\documentclass{article}

\newtheorem{theorem}{Theorem}[section]

\newtheorem{e-proposition}[theorem]{Proposition}

\newtheorem{e-definition}[theorem]{Definition\rm}
\newtheorem{remark}{\it Remark\/}
\newtheorem{example}{\it Example\/}

\title{On the flat remainder in normal forms of families of analytic planar saddles}
\author{Patrick Bonckaert and Freek Verstringe \\
Hasselt University\\ Agoralaan, gebouw D \\ B-3590 Diepenbeek, Belgium}

\begin{document}

\maketitle

\begin{abstract}
We give an explicit expression for the (finitely) flat remainder after analytic normal form reduction of a family of planar saddles of diffeomorphisms or vector fields. We distinguish between a rational or irrational ratio of the moduli of the eigenvalues at the saddle for a certain value of the parameter. 
\end{abstract}

\def\R{{\hbox{\bf R}}}
\def\N{{\hbox{\bf N}}}
\def\C{{\hbox{\bf C}}}
\def\Z{{\hbox{\bf Z}}}


\section{Introduction}
\label{sect:intro}

We consider analytic diffeomorphisms or vector fields depending on a parameter $\lambda$ and having a hyperbolic stationary point, say the origin $0$, of saddle type. For sinks or sources we only have finitely many resonances and this is a stable phenomenon when varying the parameter $\lambda$; the system can be analytically reduced to a polynomial, and this reduction depends analytically on $\lambda$ if the diffeomorphism or vector field does \cite{Bruslinskaja:1971,Gorbovitskis:2006}. This holds in any dimension.

The question of analytic reduction, for planar one-forms of saddle type, was studied in H. Dulac's memoir \cite{Dulac:1923}. The system is analytically equivalent to a normal form up to a remainder of some finite order of flatness. This form depends on the ratio of the eigenvalues at $0$. If this ratio is rational, say $-q/p$, one considers the resonant monomial $u=x_1^px_2^q$. The normal form is a function of $u$. The further elimination of the (finitely) flat remainder is of a different nature; it can, in general, only be done in a finitely smooth way \cite{Ilyashenko/Yakovenko:1991}, especially in the presence of parameters, which we shall allow in this paper. If one tries reduction by means of a formal power series, then this will generically diverge; quoting C. Rousseau: 'It is very exceptional that the change of coordinates to normal form converges'. For a geometric explanation  of this phenomenon see \cite{Rousseau:2006}. On the other hand it can be done Gevrey-one \cite{Bonckaert/DeMaesschalck:2007}. So in both cases we loose the information about the analiticity. 

We summarise that eliminating the flat remainder can remove some essential information. In this paper we give a more explicit structure for that analytic flat remainder. 

When there is no rational relation between the eigenvalues there are conditions allowing to linearize the system \cite{Siegel:1952}, but in case that there are parameters, i.e. when the eigenvalues vary, these conditions can be destroyed, so also here a flat remainder $R$ must be allowed for analytic reduction of families. 

In this paper we study planar analytic systems locally near the origin and can consider them to be sums of convergent power series on some polydisk. Although, since it is more convenient, we work in the complex, all the results remain valid if we work over $\R$, by complexifying the system. The property 'complexified system coming from a real system' is preserved at each step of the reduction below. The coefficients of the power series are functions of the parameter $\lambda$. Throughout this note they will all be defined on some sufficiently small neighbourhood $\Lambda$ of a given value $\lambda_0$.  In fact one considers these functions in a space endowed with some supremum norm. One can for example think of the space of functions that are analytic inside a polydisk and continuous on the boundary.

In the sequel we will use the multi-index conventions i.e. for $k=(k_1,k_2)\in\N^2$ and $x=(x_1,x_2)\in\C^n$ we write $x^k=x_1^{k_1} x_2^{k_2}$ and $|k|=k_1 + k_2$. 

We consider planar local analytic diffeomorphisms near $0\in\C^2$ fixing $0$. Since we will be interested in hyperbolic saddles we assume that the linear part at $0$ has eigenvalues $\mu_1,\mu_2$ with $\mu_1\ne \mu_2$. More explicitly let $D \subset \C^2$ be a polydisk at $0$ and let 
\begin{equation} \label{eq:diffeof}
f_\lambda(x_1,x_2)= (\mu_1(\lambda) x_1, \mu_2(\lambda) x_2) + \sum_{|k|=2}^\infty a_k(\lambda) x^k
\end{equation}
be absolutely convergent on $D$ where $a_k(\lambda) \in\C^2$. 

For a certain value $\lambda=\lambda_0$ of the parameter we shall distinguish between the cases that $\log |\mu_1(\lambda_0)| / \log |\mu_2(\lambda_0)|$ is rational or not. Note that this ratio can be rational, say $-q/p$, even if the eigenvalues $\mu_1(\lambda_0)$ and $\mu_2(\lambda_0)$ are nonresonant; in the real case this is equivalent with the resonance (in the usual sense) $\mu_1(\lambda_0)^{2p}=\mu_2(\lambda_0)^{-2q}$. It is standard to eliminate a finite number of nonresonant monomials by a polynomial change of variables \cite{Ilyashenko/Yakovenko:1991}. For the theorems \ref{thm:qpresonantie} and \ref{thm:irrationalratio} below we however we must consider the relation between the moduli of the eigenvalues.

\section{$q:-p$ resonance in modulus} \label{sect:qpresonance}

We say that $f_\lambda$ in (\ref{eq:diffeof}) has a $q:-p$ resonance in modulus at $\lambda=\lambda_0$ if $|\mu_1(\lambda_0)|^p=|\mu_2(\lambda_0)|^{-q}$, where $p$ and $q$ are positive integers that are relatively prime. It is known from Euclid's division algorithm and from number theory \cite{Rosen:1988} that the integer equation $qr-ps=1$ has one solution $(r,s)=(r_0,s_0)$ with $0 \le r_0 \le p$ and $0 \le s_0\le q$. We define $\alpha_0=(r_0,s_0)$ and $\alpha_1=-(r_0,s_0) + (p,q)=: (r_1,s_1)$. Observe that $qr_1-ps_1= -1$ and that also $0 \le r_1 \le p$ and $0 \le s_1 \le q$.

\begin{theorem} \label{thm:qpresonantie} Let $f_\lambda$ be as in (\ref{eq:diffeof}) and have a hyperbolic fixed point at $0$ with a $q:-p$ resonance in modulus at $\lambda=\lambda_0$. Let $N\in \N$ be given. Then there is an analytic change of variables $h_\lambda$ near $0$ such that $g_\lambda=h_\lambda^{-1} f_\lambda h_\lambda$ takes the following form, where $u=x_1^px_2^q$:
\begin{eqnarray} \label{eq:normaalvormg2}
g_\lambda(x_1,x_2) & & =  \left(x_1 [ \mu_1(\lambda) + b^0_1(u,\lambda) + u \sum_{k=1}^\infty u^{Nk} \big( x^{k\alpha_0} b^{0k}_1(u,\lambda) + x^{k\alpha_1} b^{1k}_1(u,\lambda)\big) ] , \right.
\nonumber \\
& &   x_2  \left.[ \mu_2(\lambda) + b^0_2(u,\lambda) + u \sum_{k=1}^\infty u^{Nk}\big(x^{k\alpha_0} b^{0k}_2(u,\lambda) + x^{k\alpha_1} b^{1k}_2(u,\lambda)\big) ] \right).
\end{eqnarray}
where moreover: \newline
(i) the functions $b^0_i, b^{0k}_i, b^{1k}_i$ are locally analytic in $u$, $i=1,2$; \newline
(ii) all the occurring power series expressed in the variables $x_1, x_2$ converge on some polydisk at $0$; \newline
(iii) if the family $f_\lambda$ depends analytically on $\lambda$ near $\lambda_0$ then this is also the case for $h_\lambda, g_\lambda, b^0_i, b^{0k}_i$ and $b^{1k}_i$.
\end{theorem}

\begin{example} \label{ex:p3q2} \rm For $p=3$ and $q=2$ we take $\alpha_0=(2,1)$ and $\alpha_1=(1,1)$.
\end{example}

Sketch of the \bf proof \rm. We consider the case $|\mu_1(\lambda_0)| < 1$, the other case being similar. We denote $G_0=\{(k_1(r_0,s_0) + k_2 (p,q) | k_1 \ge 1$ and $-\infty < k_2 \le (N+1)k_1\} \cap \N^2$ and $B_0=\N^2\setminus G_0$. Let us define $D_0=|\mu_1(\lambda_0)|^{1/(q(r_0 + s_0 + N(p+q)) }$
then $D_0<1$. For an element $m=(m_1,m_2) = k_1(r_0,s_0) + k_2(p,q)$ of $G_0$ we have (skipping some straightforward calculations):
\begin{eqnarray} 
|(\mu_1(\lambda_0),\mu_2(\lambda_0))^m| & = & |\mu_1(\lambda_0)|^{r_0 k_1 - \frac{p}{q} s_0 k_1} \nonumber \\
& = & D_0^{(r_0 + s_0 + (N+1)(p+q)) k_1} \nonumber \\
& \le &  D_0^{|m|} \label{eq:afschattingmum}
\end{eqnarray}
Because of (\ref{eq:afschattingmum}) we can apply the results in \cite{Bonckaert/Verstringe:2007} and infer that there exists a local analytic change of variables, which we shall denote by $\tilde h_\lambda$, with $\tilde h_\lambda(0) =0$ and tangent to the identity such that $\tilde h_{\lambda}^{-1} f_\lambda \tilde h_{\lambda} = : \tilde g_{\lambda}$ takes the form
\begin{equation} \label{tildeglambda}
\tilde g_\lambda(x_1,x_2) = (x_1(\mu_1(\lambda) + \sum_{|m|=2,m\in B_0}^\infty a_{m1}(\lambda) x^m), x_2(\mu_2(\lambda) + \sum_{|m|=2,m\in B_0}^\infty a_{m2}(\lambda) x^m)
\end{equation}
where $\tilde h_{\lambda}$, $\tilde g_{\lambda}$ depend analytically on $\lambda$ if the family $f_\lambda$ does so. All the involved power series converge absolutely on some polydisk at the origin. We can write $B_0=B_1\cup B_2$ where $B_2 = \{k_1(r_0,s_0) + k_2(p,q) | k_1 \le 0$ and $k_2 \in\Z\} \cap \N^2$ and $B_1 = \{k_1(r_0,s_0) + k_2(p,q) | k_1 \ge 1 \hbox{ and } k_2 \ge (N+1)k_1 +1\}$.
If $m\in B_1$ then $m=k_1(r_0,s_0) + k_2(p,q)$ where $k_2=(N+1)k_1 + 1 + \tilde k_2$ for some $\tilde k_2 \ge 0$. Hence, for such an $m$:
\begin{eqnarray} 
x^m & = & x^{k_1\alpha_0} u.u^{(N+1) k_1} u^{\tilde k_2}. \label{eq:xmminB1}
\end{eqnarray}
Observe that monomials in (\ref{eq:xmminB1}) are of the type of the first term of the summation $\sum_{k=1}^\infty$ in (\ref{eq:normaalvormg2}).

The next step consists in considering the inverse $\tilde g_{k}^{-1}$. Remark that $B_0$ is a subset of
\begin{equation} \label{eq:deftildeB}
\tilde B = \{(m_1,m_2) \in \N^2 | \frac{m_2}{m_1} \ge \frac{s_0 + (N+1) q}{r_0 + (N+1)p} \}.
\end{equation}
Furthermore if $m=(m_1,m_2)=k_1(r_0,s_0) + k_2(p,q)$ with $k_2 \le  Nk_1$ then 
\begin{equation} \label{eq:verhoudm2m1een} \frac{m_2}{m_1} < \frac{s_0 + (N+1) q}{r_0 + (N+1)p}.
\end{equation}
Hence for such an $m$ we have $m\in \tilde G:= \N^2 \setminus \tilde B$. So since $B_0\subset \tilde B$ we have on the formal level for the inverse $g_\lambda^{-1}$ the following expansion (see \cite{Bonckaert/Verstringe:2007}):
\begin{equation} \label{eq:formalinversetildeg}
\tilde g_{\lambda}^{-1}(x_1,x_2) = (x_1(\mu_1(\lambda)^{-1} + \sum_{|m|=2,m\in \tilde B}^\infty \hat a_{1,m}(\lambda) x^m), x_2(\mu_2(\lambda)^{-1} + \sum_{|m|=2, m\in\tilde B}^\infty \hat a_{2,m}(\lambda) x^m))
\end{equation}
Now we can reverse the role of $x_1$ and $x_2$ and observe that $|\mu_2(\lambda_0)^{-1}| < 1$. Let us denote $G_1 = \{k_1(r_1,s_1) + k_2(p,q) | k_1 \ge 1$ and $-\infty < k_2 \le (N+1)k_1\}$ and we put $D_1 = |\mu_2(\lambda_0)|^{- 1/(p(r_1+s_1+(N+1)(p+q))}$;
 observe that $D_1 < 1$. In a similar way as above we estimate that $|(\mu_1(\lambda_0)^{-1}, \mu_2(\lambda_0)^{-1})^m| \le D_1^{|m|}$ for $m\in G_1$.
Hence, using once again \cite{Bonckaert/Verstringe:2007}, there exists a local analytic change of variables $\hat h_\lambda$, with $\hat h_\lambda(0)=0$ and tangent to the identity, such that $\hat h_\lambda^{-1} \tilde g_\lambda^{-1} \hat h_\lambda = \hat g_\lambda$ takes the form
\begin{equation}\label{eq:vormhatg}
\hat g_\lambda(x_1,x_2) = (x_1(\mu_1(\lambda)^{-1} + \sum_{|m|=2,m\in \hat B}^\infty \hat b_{1,m}(\lambda) x^m), x_2(\mu_2(\lambda)^{-1} + \sum_{|m|=2,m \in \hat B}^\infty \hat b_{2,m}(\lambda) x^m))
\end{equation}
where 
\begin{equation} \label{eq:defhatB}
\hat B:=\{(m_1,m_2) \in\N^2 | \frac{s_0 + (N+1) q}{r_0 + (N+1)p} \le \frac{m_2}{m_1} \le \frac{s_1+(N+1)q}{r_1+(N+1)p} \}.
\end{equation} 
Moreover also $\hat g_\lambda^{-1}$ takes the form
\begin{equation}\label{eq:vormhatginvers}
\hat g_\lambda^{-1}(x_1,x_2) = (x_1(\mu_1(\lambda)^{-1} + \sum_{|m|=2,m\in \hat B}^\infty b_{1,m}(\lambda) x^m), x_2(\mu_2(\lambda)^{-1} + \sum_{|m|=2,m \in \hat B}^\infty b_{2,m}(\lambda) x^m)).
\end{equation}
Obviously, $\hat h_\lambda^{-1}\tilde h_\lambda^{-1} f_\lambda \tilde h_\lambda \hat h_\lambda = \hat g_\lambda^{-1}$, which means that $f_\lambda$ is conjugate to $\hat g_\lambda^{-1}$ by $\tilde h_\lambda \circ \hat h_\lambda$. We put $g_\lambda=\hat g_\lambda^{-1}$. It is not hard to see that $g_\lambda$ has the aimed form (\ref{eq:normaalvormg2}), because of the choice of $(r_0,s_0)$. 

\section{Nonresonance}

Let $f_\lambda$ be as in (\ref{eq:diffeof}) with 
\begin{equation} \label{eq:defR}
\frac{\log |\mu_1(\lambda_0)|}{\log |\mu_2(\lambda_0)|} = -R
\end{equation}
with $R>0$. If $R$ is irrational we consider the continued fraction expansion 
\begin{equation} \label{eq:kettingbreukR}
R=[a_1\  a_2\  \dots a_n \dots]
\end{equation}
of $R$ (see for instance \cite{Olds:1963}, also for the statements (\ref{eq:benaderingdoorconvergenten}) and (\ref{eq:qnpngelijkeen}) made below). Let 
\begin{equation}\label{eq:convergenten}
\frac{q_n}{p_n} = [a_1\ a_2 \dots a_n]
\end{equation}
denote the convergents of $R$. One has:
\begin{equation} \label{eq:benaderingdoorconvergenten}
\frac{q_1}{p_1} < \frac{q_3}{p_3} < \dots < \frac{q_{2k+1}}{p_{2k+1}} < \dots < R < \dots < \frac{q_{2k+2}}{p_{2k+2}} < \dots < \frac{q_2}{p_2}.
\end{equation}
The sequences $(p_n)_{n\in\N}$ and $(q_n)_{n\in\N}$ tend to infinity at least as fast as the Fibonacci numbers. Our intention is to take $p_n$ and $q_n$ 'as large as desired'. We will crucially make use of the known fact that
\begin{equation}\label{eq:qnpngelijkeen}
q_{n}p_{n+1} - q_{n+1}p_n=(-1)^{n+1}.
\end{equation}

\begin{theorem} \label{thm:irrationalratio} Let $f_\lambda$ and $R$ be as in (\ref{eq:diffeof}) respectively (\ref{eq:defR}), and suppose that $0$ is a hyperbolic fixed point of $f_\lambda$. Let $k\in\N$ be given and denote $q=q_{2k+1}$ and $p=p_{2k+1}$, using the notation from (\ref{eq:kettingbreukR}) and (\ref{eq:convergenten}). Denote $\tilde q=q_{2k+2}$ and $\tilde p=p_{2k+2}$. Then there  is an analytic change of variables $h_\lambda$ near $0$ such that $g_\lambda=h_\lambda^{-1}f_\lambda h_\lambda$ takes the following form, where $u=x_1^px_2^q$ and $\tilde u=x_1^{\tilde p}x_2^{\tilde q}$:
\begin{eqnarray} 
g_\lambda(x_1,x_2) & = & (x_1[\mu_1(\lambda) + ub^1_1(u,\tilde u,\lambda) + \tilde ub^2_1(u,\tilde u,\lambda)] , \nonumber \\
& & x_2[\mu_2(\lambda) + u b^1_2(u,\tilde u,\lambda) + \tilde u b^2_2(u,\tilde u,\lambda) ] ) \label{eq:normaalvormirrational}
\end{eqnarray}
where moreover: \newline
(i) the functions $b^i_j$ are locally analytic near $(u,\tilde u)=(0,0)$, $i,j=1,2$; \newline
(ii) if the family $f_\lambda$ depends analytically on $\lambda$ near $\lambda_0$ then this is also the case for $h_\lambda$, $g_\lambda$, $b^i_j$, $i,j=1,2$.
\end{theorem}

The \bf proof \rm of this theorem follows similar lines as in section \ref{sect:qpresonance}; here we take $D_0=|\mu_1(\lambda_0)|^{(Rp-q)/(R(p+q))}$ and $D_1=|\mu_2(\lambda_0)^{-1}|^{(\tilde q - R \tilde p)/(\tilde p + \tilde q)}$. 

\begin{example} \rm
Let $R=\phi = (1+\sqrt{5})/2$ be the golden ratio and $(F_n)_{n\in\N}$ the Fibonacci sequence. Then we take $u=x_1^{F_{2k+1}}x_2^{F_{2k+2}}$ and $\tilde u=u x_1^{F_{2k}}x_2^{F_{2k+1}}$ in expression (\ref{eq:normaalvormirrational}).
\end{example}

\section{Final remarks}

\begin{remark} In theorems  \ref{thm:qpresonantie} and \ref{thm:irrationalratio} we have that the axes $\{x_1=0\}$ and $\{x_2=0\}$ are invariant for $g_\lambda$. This implies the classical fact that a hyperbolic analytic (family-) germ has analytic stable and unstable manifolds. Moreover, from formulas (\ref{eq:normaalvormg2}) and  (\ref{eq:normaalvormirrational}) we see that the flat remainder of these normal forms $g_\lambda$ is as (finitely) flat as desired along both $\{x_1=0\}$ and $\{x_2=0\}$. In fact this holds in any dimension.
\end{remark}

\begin{remark} Using the results in \cite{Bonckaert/DeMaesschalck:2007} we can formulate analogous theorems for saddle singularities of planar vector fields.
\end{remark}





\end{document}